\documentclass[a4paper,12pt,dvips]{article} 
\usepackage{amssymb}
\usepackage{amscd}
\usepackage{amsmath} 
\usepackage{graphicx}
\usepackage{latexsym} 
\usepackage{enumerate}

\usepackage{theorem}
\theoremstyle{plain}
\theorembodyfont{\itshape}
\newtheorem{theorem}{Theorem}[section]
\newtheorem{proposition}[theorem]{Proposition}
\newtheorem{lemma}[theorem]{Lemma}
\newtheorem{corollary}[theorem]{Corollary} 
\theorembodyfont{\rmfamily}

\title{\bf Cubic Harmonics and Bernoulli Numbers} 
\author{Katsunori Iwasaki\thanks{Department of Mathematics, 
Hokkaido University, Kita 10, Nishi 8, Kita-ku, Sapporo 060-0810 Japan. 
E-mail: {\tt iwasaki@math.sci.hokudai.ac.jp}}}
\date{October 25, 2011} 
\begin{document}
\maketitle
\begin{abstract}
The functions satisfying the mean value property for an 
$n$-dimensional cube are determined explicitly. 
This problem is related to invariant theory for a finite 
reflection group, especially to a system of invariant 
differential equations.  
Solving this problem is reduced to showing that a certain set 
of invariant polynomials forms an invariant basis. 
After establishing a certain summation formula over Young diagrams, 
the latter problem is settled by considering  a recursion formula 
involving Bernoulli numbers. \\[3mm]
\end{abstract} 
\section{Introduction} \label{sec:intro} 
Let $P$ be an $n$-dimensional polytope in $\mathbb{R}^n$. 
For $k = 0,\dots,n$, let $P(k)$ be the $k$-dimensional skeleton of $P$.  
A continuous function $f : \mathbb{R}^n \to \mathbb{R}$ is said to be 
{\it $P(k)$-harmonic} if it satisfies 
\begin{equation} \label{eqn:mvp}
f(x) = \dfrac{1}{|P(k)|} \int_{P(k)} f(x + r y) \, d\mu_k(y) 
\end{equation}
for any $x \in \mathbb{R}^n$ and $r > 0$, where $\mu_k$ is the $k$-dimensional 
Euclidean measure on $P(k)$ and $|P(k)| := \mu_k(P(k))$ is the 
$k$-dimensional Euclidean volume of $P(k)$. 
This is an extension to a polytope of the classical notion of 
harmonic functions characterized by the mean value 
property for the $(n-1)$-dimensional sphere $S^{n-1}$.  
Let $\mathcal{H}_{P(k)}$ denote the set of all $P(k)$-harmonic functions on 
$\mathbb{R}^n$.  
A general result in Iwasaki \cite{Iwasaki} states that for any 
polytope $P$ and any $k = 0,\dots, n$, the set $\mathcal{H}_{P(k)}$ is 
a finite-dimensional linear space of polynomials. 
Note that $\mathcal{H}_{P(k)}$ carries the structure of an $\mathbb{R}[\partial]$-module, 
because equation (\ref{eqn:mvp}) is stable under partial differentiations  
$\partial = (\partial_1, \dots, \partial_n)$, where   
$\partial_i := \partial/\partial x_i$ is the $i$-th partial 
differential operator.    
\par
It is an interesting problem to determine the space $\mathcal{H}_{P(k)}$ 
explicitly when $P$ is a regular convex polytope with center at 
the origin in $\mathbb{R}^n$.  
This problem is already settled unless $P$ is an 
$n$-dimensional cube.  
As for the cube case, however, although the vertex problem $(k = 0)$ was 
solved by Flatto \cite{Flatto2} and Haeuslein \cite{Haeuslein} as 
early as 1970, the higher skeleton problem $(k = 1, \dots, n)$ has 
been open up to now. 
The aim of this article is to give a complete solution to 
this problem. 
A characteristic feature of our work is that it provides a  
simultaneous resolution for all skeletons, which reveals a natural 
structure of this problem from the viewpoint of combinatorial analysis. 
Here we refer to Iwasaki \cite{Iwasaki3} for a general review of 
the topic discussed in this article.      
\par
Let $C$ be an $n$-dimensional cube with center at the origin in 
the Euclidean space $\mathbb{R}^n$ endowed with the standard orthonormal 
coordinates $x = (x_1,\dots,x_n)$.  
After a scale change and a rotation one may assume that the 
vertices of $C$ are at $(\pm1,\dots,\pm1)$.  
The symmetry group of $C$ is a finite reflection group of type $B_n$, 
which is the semi-direct product $W_n := S_n \ltimes \{\pm 1\}^n$ of  
the group $\{\pm 1\}^n = \{\varepsilon = (\varepsilon_1,\dots,\varepsilon_n)\,:\, \varepsilon_i = \pm1\}$ 
of $n$-tuple signs acting on $\mathbb{R}^n$ by sign changes of $x_1, \dots, x_n$, 
with the symmetric group $S_n$ acting on $\mathbb{R}^n$ by 
permuting $x_1, \dots, x_n$. 
The order of $W_n$ is $2^n \cdot n!$ and the fundamental alternating 
polynomial of $W_n$ is given by 
\[
\varDelta(x) = x_1 \cdots x_n \prod_{i < j} (x_i^2 - x_j^2). 
\] 
The first main result of this article is then stated as follows.  
\begin{theorem} \label{thm:harmonic} 
Let $C$ be an $n$-dimensional cube centered at the origin in 
$\mathbb{R}^n$. 
For any $k = 0, \dots, n$, the linear space $\mathcal{H}_{C(k)}$ is of  
$2^n \cdot n!$-dimensions and as an $\mathbb{R}[\partial]$-module 
$\mathcal{H}_{P(k)}$ is generated by the fundamental alternating 
polynomial $\varDelta(x)$ of the reflection group $W_n$.   
\end{theorem}
\par
For an arbitrary polytope $P$ Iwasaki \cite{Iwasaki} introduced an 
infinite sequence of homogeneous polynomials $\tau_m^{(k)}(x)$ of 
degrees $m = 1,2,3,\dots$ in terms of some combinatorial data about 
$P(k)$ and characterized $\mathcal{H}_{P(k)}$ as 
the solution space to the system of partial differential equations
\[
\tau_m^{(k)}(\partial) f = 0 \qquad (m = 1,2,3, \dots).  
\]
From the way in which they are defined, the polynomials $\tau_m^{(k)}(x)$ 
are invariant under the symmetry group $G$ of $P$. 
This observation connects our problem to the theory of $G$-harmonic 
functions due to Steinberg \cite{Steinberg}. 
A $C^{\infty}$-function $f : \mathbb{R}^n \to \mathbb{R}$ is said to be {\it $G$-harmonic} 
if it satisfies $\varphi(\partial) f = 0$ for any  
$G$-invariant polynomial $\varphi(x)$ without constant term.  
Let $\mathcal{H}_G$ denote the set of all $G$-harmonic functions. 
There is always the inclusion $\mathcal{H}_G \subset \mathcal{H}_{P(k)}$, and 
if $\{\tau_m^{(k)}(x)\}_{m=1}^{\infty}$ happens to generate the 
ring of $G$-invariant polynomials, then there occurs the coincidence 
$\mathcal{H}_G = \mathcal{H}_{P(k)}$. 
Steinberg \cite{Steinberg} made an explicit determination of $\mathcal{H}_G$ 
when $G$ is a finite reflection group: $\dim \mathcal{H}_G = |G|$ and 
as an $\mathbb{R}[\partial]$-module $\mathcal{H}_G$ is generated by the fundamental 
alternating polynomial of $G$. 
If $P$ is a regular polytope, then $G$ is a finite 
reflection group and we are done if we are able to show that the 
sequence $\{\tau_m^{(k)}(x)\}_{m=1}^{\infty}$ actually generates  
the $G$-invariant ring.           
\par
For the $k$-skeleton $C(k)$ of the $n$-cube $C$ the polynomials 
$\tau_m^{(k)}(x)$ are constructed as follows.     
First recall that the $m$-th complete symmetric polynomial of 
$j$ variables is defined by   
\begin{equation} \label{eqn:H}
H_m^{(j)}(t_1,\dots,t_j) := 
\sum_{(m_1, \dots, m_j) \in \mathcal{P}_j(m)} t_1^{m_1} \cdots t_j^{m_j},  
\end{equation}
where $\mathcal{P}_j(m)$ is the set of all ordered partitions  
$(m_1,\dots,m_j)$ of $m$ by $j$ nonnegative integers.  
Note that $H_m^{(j+1)}(t_1, \dots, t_j, 0) = H_m^{(j)}(t_1, \dots, t_j)$.   
For each $k = 0,\dots,n$, we next define    
\begin{equation} \label{eqn:h}
h_m^{(k)}(x) := 
H_m^{(k+1)}(x_1+\cdots+x_n, \,\, x_2+\cdots+x_n, \,\, 
\dots, \,\, x_{k+1}+\cdots+x_n).     
\end{equation}
Note that when $k = n$ the term $x_{k+1}+\cdots+x_n$ is null and 
thus $h_m^{(n)}(x) = h_m^{(n-1)}(x)$. 
\par
\begin{figure}[t]
\begin{center}
\unitlength 0.1in
\begin{picture}( 37.3600, 35.7800)(  2.1000,-38.9400)
%
\special{pn 20}%
\special{pa 1322 870}%
\special{pa 362 1510}%
\special{fp}%
\special{pa 362 1510}%
\special{pa 362 1510}%
\special{fp}%
%
\special{pn 20}%
\special{pa 3562 878}%
\special{pa 2602 1518}%
\special{fp}%
\special{pa 2602 1518}%
\special{pa 2602 1518}%
\special{fp}%
%
\special{pn 20}%
\special{pa 362 1510}%
\special{pa 2602 1510}%
\special{fp}%
%
\special{pn 20}%
\special{pa 362 1510}%
\special{pa 362 3750}%
\special{fp}%
%
\special{pn 20}%
\special{pa 2610 1510}%
\special{pa 2610 3750}%
\special{fp}%
%
\special{pn 20}%
\special{pa 362 3750}%
\special{pa 2602 3750}%
\special{fp}%
%
\special{pn 20}%
\special{pa 3562 3118}%
\special{pa 2602 3758}%
\special{fp}%
\special{pa 2602 3758}%
\special{pa 2602 3758}%
\special{fp}%
%
\special{pn 20}%
\special{pa 1322 3102}%
\special{pa 362 3742}%
\special{dt 0.054}%
\special{pa 362 3742}%
\special{pa 362 3742}%
\special{dt 0.054}%
%
\special{pn 20}%
\special{sh 1.000}%
\special{ia 2442 870 32 32  0.0000000 6.2831853}%
\special{ar 2442 870 32 32  0.0000000 0.3870968}%
\special{ar 2442 870 32 32  1.5483871 1.9354839}%
\special{ar 2442 870 32 32  3.0967742 3.4838710}%
\special{ar 2442 870 32 32  4.6451613 5.0322581}%
\special{ar 2442 870 32 32  6.1935484 6.2832853}%
%
\special{pn 20}%
\special{sh 1.000}%
\special{ia 1962 1190 32 32  0.0000000 6.2831853}%
\special{ar 1962 1190 32 32  0.0000000 0.3870968}%
\special{ar 1962 1190 32 32  1.5483871 1.9354839}%
\special{ar 1962 1190 32 32  3.0967742 3.4838710}%
\special{ar 1962 1190 32 32  4.6451613 5.0322581}%
\special{ar 1962 1190 32 32  6.1935484 6.2832853}%
%
\special{pn 20}%
\special{sh 1.000}%
\special{ia 3562 870 32 32  0.0000000 6.2831853}%
\special{ar 3562 870 32 32  0.0000000 0.3870968}%
\special{ar 3562 870 32 32  1.5483871 1.9354839}%
\special{ar 3562 870 32 32  3.0967742 3.4838710}%
\special{ar 3562 870 32 32  4.6451613 5.0322581}%
\special{ar 3562 870 32 32  6.1935484 6.2832853}%
%
\special{pn 13}%
\special{pa 1962 1182}%
\special{pa 2450 870}%
\special{fp}%
%
\special{pn 20}%
\special{sh 1.000}%
\special{ia 1962 2310 32 32  0.0000000 6.2831853}%
\special{ar 1962 2310 32 32  0.0000000 0.3870968}%
\special{ar 1962 2310 32 32  1.5483871 1.9354839}%
\special{ar 1962 2310 32 32  3.0967742 3.4838710}%
\special{ar 1962 2310 32 32  4.6451613 5.0322581}%
\special{ar 1962 2310 32 32  6.1935484 6.2832853}%
%
\special{pn 13}%
\special{pa 1962 1198}%
\special{pa 3562 862}%
\special{fp}%
\put(20.1800,-24.8600){\makebox(0,0)[lb]{O}}%
\put(17.7000,-12.4000){\makebox(0,0)[lb]{F}}%
\put(25.2000,-7.8000){\makebox(0,0)[lb]{E}}%
\put(36.6000,-8.5000){\makebox(0,0)[lb]{V}}%
%
\special{pn 13}%
\special{sh 1.000}%
\special{ar 3082 2310 20 20  0.0000000 6.2831853}%
%
\special{pn 13}%
\special{sh 1.000}%
\special{ar 1962 3430 20 20  0.0000000 6.2831853}%
%
\special{pn 13}%
\special{sh 1.000}%
\special{ar 850 2310 20 20  0.0000000 6.2831853}%
\put(39.4600,-23.6600){\makebox(0,0)[lb]{$x_1$}}%
\put(31.4600,-15.5000){\makebox(0,0)[lb]{$x_2$}}%
\put(18.7400,-4.8600){\makebox(0,0)[lb]{$x_3$}}%
%
\special{pn 8}%
\special{pa 1962 1566}%
\special{pa 1962 3654}%
\special{fp}%
%
\special{pn 8}%
\special{pa 1962 1446}%
\special{pa 1962 566}%
\special{fp}%
\special{sh 1}%
\special{pa 1962 566}%
\special{pa 1942 634}%
\special{pa 1962 620}%
\special{pa 1982 634}%
\special{pa 1962 566}%
\special{fp}%
%
\special{pn 8}%
\special{pa 1970 3894}%
\special{pa 1970 3798}%
\special{fp}%
%
\special{pn 8}%
\special{pa 434 2310}%
\special{pa 2546 2310}%
\special{fp}%
%
\special{pn 8}%
\special{pa 2666 2302}%
\special{pa 3882 2302}%
\special{fp}%
\special{sh 1}%
\special{pa 3882 2302}%
\special{pa 3816 2282}%
\special{pa 3830 2302}%
\special{pa 3816 2322}%
\special{pa 3882 2302}%
\special{fp}%
%
\special{pn 8}%
\special{pa 210 2310}%
\special{pa 306 2310}%
\special{fp}%
%
\special{pn 8}%
\special{pa 1402 2678}%
\special{pa 2554 1910}%
\special{fp}%
%
\special{pn 8}%
\special{pa 2658 1846}%
\special{pa 3122 1534}%
\special{fp}%
\special{sh 1}%
\special{pa 3122 1534}%
\special{pa 3056 1556}%
\special{pa 3078 1564}%
\special{pa 3078 1588}%
\special{pa 3122 1534}%
\special{fp}%
%
\special{pn 20}%
\special{pa 3562 2254}%
\special{pa 3562 870}%
\special{fp}%
%
\special{pn 20}%
\special{pa 3562 3102}%
\special{pa 3562 2334}%
\special{fp}%
%
\special{pn 20}%
\special{pa 1322 870}%
\special{pa 1930 870}%
\special{fp}%
%
\special{pn 20}%
\special{pa 2002 862}%
\special{pa 3554 862}%
\special{fp}%
%
\special{pn 20}%
\special{pa 1322 870}%
\special{pa 1322 1462}%
\special{dt 0.054}%
%
\special{pn 20}%
\special{pa 1322 1558}%
\special{pa 1322 2262}%
\special{dt 0.054}%
%
\special{pn 20}%
\special{pa 1322 2358}%
\special{pa 1322 3110}%
\special{dt 0.054}%
%
\special{pn 20}%
\special{pa 1322 3110}%
\special{pa 1922 3110}%
\special{dt 0.054}%
%
\special{pn 20}%
\special{pa 2002 3110}%
\special{pa 2562 3110}%
\special{dt 0.054}%
%
\special{pn 20}%
\special{pa 2666 3110}%
\special{pa 3562 3110}%
\special{dt 0.054}%
%
\special{pn 13}%
\special{pa 1962 2310}%
\special{pa 2562 1766}%
\special{dt 0.045}%
%
\special{pn 13}%
\special{pa 2650 1686}%
\special{pa 3570 862}%
\special{dt 0.045}%
%
\special{pn 13}%
\special{pa 1962 2310}%
\special{pa 2210 1550}%
\special{dt 0.045}%
%
\special{pn 13}%
\special{pa 2250 1454}%
\special{pa 2354 1150}%
\special{dt 0.045}%
%
\special{pn 13}%
\special{pa 2370 1070}%
\special{pa 2442 870}%
\special{dt 0.045}%
\end{picture}%
\end{center}
\caption{The cube in three dimensions with a 
flag $\mathrm{V} \prec \mathrm{E} \prec \mathrm{F} \prec \mathrm{O}$}
\label{fig:cube}
\end{figure}
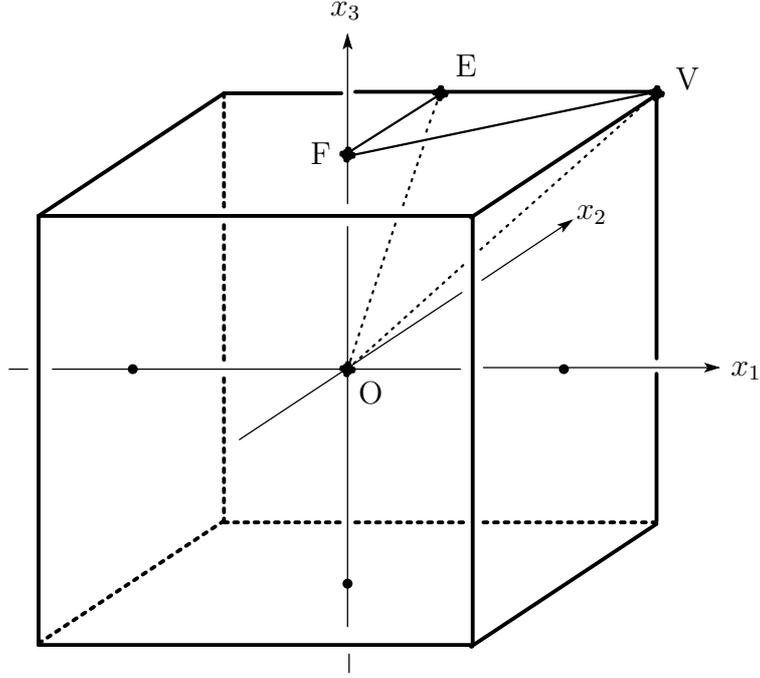

 
For example, when $n = 3$ these polynomials are given by 
\[
h_m^{(k)}(x) = \left\{\begin{array}{ll}
H_m^{(1)}(\mathrm{V} \cdot x) \qquad & (k = 0), \\[2mm]
H_m^{(2)}(\mathrm{V} \cdot x, \mathrm{E} \cdot x) \qquad & (k = 1), \\[2mm]
H_m^{(3)}(\mathrm{V} \cdot x, \, \mathrm{E} \cdot x, \, \mathrm{F} \cdot x) \qquad & (k = 2, 3),  
\end{array} 
\right.
\]
where $\mathrm{V} = (1,1,1)$ is a vertex, $\mathrm{E} = (0,1,1)$ is the midpoint of 
an edge and $\mathrm{F} = (0,0,1)$ is the center of a face of the $3$-cube $C$  
(see Figure \ref{fig:cube}) and $\mathrm{V}\cdot x$ stands for the 
inner product of $\mathrm{V}$ and $x$ regarded as space vectors.  
Identify $\mathrm{E}$ and $\mathrm{F}$ with the edge and the face on which they lie.    
Similarly the origin $\mathrm{O}$, i.e., the center of the cube $C$ is 
identified with the unique $3$-cell, i.e., the cube itself.    
Then we have a flag $\mathrm{V} \prec \mathrm{E} \prec \mathrm{F} \prec \mathrm{O}$, 
where $* \prec **$ indicates that $*$ is a face of $**$. 
The simplex $\mathrm{V} \mathrm{E} \mathrm{F} \mathrm{O}$ is a fundamental domain of 
the symmetry group $W_3$. 
There is a bijection between the elements of $W_3$ and the 
flags of $C$.  
These pictures carry over in $n$ dimensions.     
\par
Finally $\tau_m^{(k)}(x)$ is defined to be the $W_n$-symmetrization of 
$h_m^{(k)}(x)$, that is, the average: 
\begin{equation} \label{eqn:tau} 
\tau_m^{(k)}(x) := \dfrac{1}{2^n \cdot n!} \sum_{\sigma \in W_n} 
h_m^{(k)}(\sigma x) \qquad (k = 0,\dots,n). 
\end{equation}
In other words $\tau_m^{(k)}(x)$ is the average of $h_m^{(k)}(x)$ over 
all flags of $C$. 
Note that $\tau_m^{(n)}(x) = \tau_m^{(n-1)}(x)$ since 
$h_m^{(n)}(x) = h_m^{(n-1)}(x)$ as mentioned earlier. 
It is immediate from definition (\ref{eqn:tau}) that 
$\tau_m^{(k)}(x)$ is a homogeneous $W_n$-invariant of degree $m$. 
Recall that the degrees of $W_n$ are $2, 4, \dots, 2n$, 
which are all even (see e.g. Humphreys \cite{Humphreys}).   
So the invariant polynomial $\tau_m^{(k)}(x)$ vanishes identically for 
every $m$ odd. 
The second main result of this article is then stated as follows. 
\begin{theorem} \label{thm:basis} 
For any $k = 0, \dots, n$, the polynomials $\tau_2^{(k)}(x), 
\tau_4^{(k)}(x), \dots, \tau_{2n}^{(k)}(x)$ form an invariant 
basis of the reflection group $W_n$.   
\end{theorem} 
\par
For the proof of this we recall that 
$e_2(x), \dots, e_{2n}(x)$ form an invariant basis of $W_n$, where   
\[
e_{2m}(x) := \sum_{1 \le i_1 < \cdots < i_m \le n} 
x_{i_1}^2 \cdots x_{i_m}^2  \qquad (m = 1, \dots, n) 
\]
is the $m$-th elementary symmetric polynomial of $x_1^2, \dots, x_n^2$. 
Since $\tau_{2m}^{(k)}(x)$ is a homogeneous $W_n$-invariant of degree 
$2m$, there exist a unique constant $c_{n,m}^{(k)}$ and a unique weighted 
homogeneous polynomial $P_{n,m}^{(k)}(t_1, \dots, t_{m-1})$ of degree 
$2m$ with $t_i$ being of weight $2i$ such that 
\begin{equation} \label{eqn:tau-s}
\tau_{2m}^{(k)}(x) = c_{n,m}^{(k)} \, e_{2m}(x) + 
P_{n,m}^{(k)}(e_2(x), \dots, e_{2m-2}(x)) 
\qquad (m = 1,\dots,n).     
\end{equation}
Here we employ the notation $c_{n,m}^{(k)}$ and $P_{n,m}^{(k)}$ 
to emphasize the dependence upon $n$. 
Note that $c_{n,m}^{(n)} = c_{n,m}^{(n-1)}$ since 
$\tau_{2m}^{(n)}(x) = \tau_{2m}^{(n-1)}(x)$. 
If we are able to show that the coefficient $c_{n,m}^{(k)}$ does not 
vanish for any $m = 1,\dots,n$, then we can invert equations 
(\ref{eqn:tau-s}) to express $e_2(x), \dots, e_{2m}(x)$ as polynomials 
of $\tau_{2}^{(k)}(x), \dots, \tau_{2m}^{(k)}(x)$. 
From this Theorem \ref{thm:basis} follows immediately. 
So it is important to develop a method to calculate 
$c_{n,m}^{(k)}$ or at least to show that it does not vanish.    
\par 
It turns out that the coefficients $c_{n,m}^{(k)}$ exhibit a beautiful 
combinatorial structure upon introducing the generating polynomials 
\begin{equation} \label{eqn:Gnm}
G_{n,m}(t) := \sum_{k=0}^n \dfrac{n! \, c_{n,m}^{(k)}}{(n-k)! \, (2m+k)!} \,  
t^{n-k}.    
\end{equation}
The third main result of this article is concerned with the 
structure of these polynomials.    
\begin{theorem} \label{thm:generating} 
The polynomials $G_{n,m}(t)$ are tied to $G_m(t) := G_{m,m}(t)$ 
by a simple relation 
\begin{equation} \label{eqn:Gnm-Gmm}
G_{n,m}(t) = (t+1)^{n-m} G_m(t) \qquad (n \ge m \ge 1). 
\end{equation}
On the other hand the polynomials $G_m(t)$ admit a 
generating series representation   
\begin{equation} \label{eqn:Gmm}
\sum_{m=1}^{\infty} 
(-1)^{m-1} G_m(t) \left(\dfrac{z^2}{t+1}\right)^m 
= \dfrac{z \coth z + t z^2 -1}{2(tz \coth z + 1)}. 
\end{equation}
\end{theorem}
\par
Equation (\ref{eqn:Gmm}) readily leads to a recursion formula for 
$G_m(t)$ involving the Bernoulli numbers $B_m$.  
There are several conventions for defining Bernoulli numbers, but 
the most useful one in our context is through the Maclaurin series 
expansion 
\[
\dfrac{z}{e^z-1} = 1 - \dfrac{z}{2} + \sum_{m=1}^{\infty} (-1)^{m-1} 
\dfrac{B_m}{(2m)!} z^{2m},   
\]
or equivalently through the formula 
\begin{equation} \label{eqn:zcothz}
z \coth z = 1 + 2 \sum_{m=1}^{\infty} (-1)^{m-1} b_m z^{2m}, 
\qquad b_m :=  \dfrac{2^{2m-1}}{(2m)!} B_m.    
\end{equation}
Multiplying formula (\ref{eqn:Gmm}) by $2(tz \coth z +1)$, expanding 
the resulting equation into a power series of $z^2$, and comparing 
the $m$-th coefficients of both sides, we obtain the following.  
\begin{corollary} \label{cor:recursion} 
The polynomials $G_m(t)$ satisfy a recursion formula 
\begin{equation} \label{eqn:recursion} 
G_m(t) = b_m (t+1)^{m-1} + 2t \sum_{i=1}^{m-1} b_{m-i} 
(t+1)^{m-i-1} G_i(t), \qquad 
G_1(t) = \dfrac{t}{2} + \dfrac{1}{6}. 
\end{equation}
\end{corollary}
\par
A polynomial of degree $m$ is said to be {\it positive} if its 
coefficients up to degree $m$ are all positive. 
Note that the product of a positive polynomial of degree $i$ 
and a positive polynomial of degree $j$ is a positive polynomial 
of degree $i+j$. 
With definition (\ref{eqn:zcothz}) we have    
\[
b_m = \dfrac{\zeta(2m)}{\pi^{2m}} = \dfrac{1}{\pi^{2m}} 
\sum_{j=1}^{\infty} \dfrac{1}{j^{2m}} > 0 \qquad (m = 1,2,3,\dots). 
\]
Therefore recursion formula (\ref{eqn:recursion}) inductively 
implies that $G_m(t)$ is a positive polynomial of degree $m$. 
Formula (\ref{eqn:Gnm-Gmm}) then tells us that $G_{n,m}(t)$ is 
a positive polynomial of degree $n$.  
Finally formula (\ref{eqn:Gnm}) concludes that the coefficient 
$c_{m,n}^{(k)}$ is positive for any $n \ge m \ge 1$ and 
$k = 0, \dots, n$.  
This establishes Theorem \ref{thm:basis}. 
The logical structure of our main results is this: 
\[
\begin{CD}
\mbox{Theorem \ref{thm:generating}} @>>> 
\mbox{Corollary \ref{cor:recursion}} @>>>
\mbox{Theorem \ref{thm:basis}} @>>> 
\mbox{Theorem \ref{thm:harmonic}}. 
\end{CD}
\]
Thus the main body of this article is exclusively devoted to 
establishing Theorem \ref{thm:generating}.  
\par
The plan of this article is as follows. 
In Section \ref{sec:matrix} we represent the coefficient  
$c_{n,m}^{(k)}$ in terms of a sum over some matrices 
(see Proposition \ref{prop:c3}). 
In Section \ref{sec:young} this representation is recast to a 
summation formula over some Young diagrams (see 
Proposition \ref{prop:c5}).  
After these preliminaries, Theorem \ref{thm:generating} 
and Corollary \ref{cor:recursion} are established in 
Section \ref{sec:generating}, where some amplifications 
of these results and a summary on polyhedral harmonics 
for regular convex polytopes are also included.     
\section{Matrix Representation} \label{sec:matrix} 
We derive a representation of the coefficient $c_{n,m}^{(k)}$ as 
the sum of some quantities depending on a certain class of matrices. 
The main result of this section is given in Proposition \ref{prop:c3}.    
Various representations in this section involve those matrices as 
in Figure \ref{fig:upper}, namely, $A = (a_{ij})$ with 
$a_{ij} = 0$ for any $i > j$. 
Such a matrix is referred to as an {\it upper quadrilateral} matrix. 
\begin{figure} 
\[
A = \begin{pmatrix}
* & * & * & * & * & * & * \\
  & * & * & * & * & * & * \\
  &   & * & * & * & * & * \\
  &   &   & * & * & * & * \\
  &   &   &   & * & * & *  
\end{pmatrix}
\]
\caption{An upper quadrilateral matrix}
\label{fig:upper}
\end{figure}
Note that it becomes an upper triangular matrix if its vertical size 
is larger than or equal to its horizontal size.  
Throughout this article we use the following notation. 
For a matrix $M = (m_{ij})$ of nonnegative integers, whether upper 
quadrilateral or not, or even for a row or column vector,    
\[
M! := \prod_{i,j} m_{ij}!.  
\]
For a row vector $\vec{v} = (v_1, \dots, v_n)$ of nonnegative integers  
we put $x^{\vec{v}} := x_1^{v_1} \cdots x_n^{v_n}$. 
Moreover,    
\[
\vec{e} := (\overbrace{1, \, \dots, \,1}^{k+1}), \qquad 
\mathbf{1} := \left. \begin{pmatrix} 1 \\ \vdots \\ 1\end{pmatrix} 
\right\}\!{\scriptstyle n}. 
\]
\begin{lemma} \label{lem:h} The polynomial in $(\ref{eqn:h})$ is 
expressed as  
\begin{equation} \label{eqn:h-M}
h_m^{(k)}(x) = \sum_{A} \dfrac{(A \mathbf{1})!}{A!} x^{\vec{e} A}, 
\end{equation}
where the sum is taken over all $(k+1) \times n$ quadrilateral 
matrices $A$ of nonnegative integers whose entries 
sum up to $m$.    
\end{lemma}
{\it Proof}. In view of definitions (\ref{eqn:H}) and (\ref{eqn:h}),  
the multi-nomial theorem yields
\[
\begin{array}{rcl}
h_m^{(k)}(x) &=& \displaystyle \sum_{m_1+\cdots+m_{k+1}=m} 
\displaystyle \prod_{i=1}^{k+1}(x_i+\cdots+x_n)^{m_i} \\[3mm]
&=& \displaystyle \sum_{m_1+\cdots+m_{k+1}=m} 
\prod_{i=1}^{k+1} \left( 
\sum_{a_{ii}+\cdots+a_{in}=m_i} 
\dfrac{m_i!}{\prod_{j=i}^n a_{ij}!} \prod_{j=i}^n x_j^{a_{ij}} 
\right) \\[8mm]
&=& \displaystyle \sum_{\sum a_{ij} = m} 
\dfrac{\prod_{i=1}^{k+1} (a_{ii}+\cdots+a_{in})!}{\prod_{i\le j} a_{ij}!} 
\prod_{i \le j} x_j^{a_{ij}} \\[8mm]
&=& \displaystyle \sum_{\sum a_{ij} = m} 
\dfrac{\prod_{i=1}^{k+1} (a_{ii}+\cdots+a_{in})!}{\prod_{i\le j} a_{ij}!} 
\prod_{j=1}^n x_j^{\sum_{i=1}^{\min\{j,k+1 \}} a_{ij}},   
\end{array}
\]
where $a_{ij}$ is defined for $1 \le i \le j \le n$, $i \le k+1$. 
Putting $a_{ij} = 0$ for $k+1 \ge i > j \ge 1$ makes 
$A =(a_{ij})$ an upper quadrilateral matrix.    
It is obvious that the entries of $A$ sum up to $m$.  \hfill $\Box$  
\par\medskip
Consider the $\{\pm 1\}^n$-symmetrization of $h_m^{(k)}(x)$, that is, 
the average:    
\begin{equation} \label{eqn:g} 
g_m^{(k)}(x) := \dfrac{1}{2^n}
\sum_{\varepsilon \in \{\pm 1\}^n} h_m^{(k)}(\varepsilon_1 x_1, \dots, \varepsilon_n x_n), 
\qquad 
\varepsilon = (\varepsilon_1, \dots, \varepsilon_n) \in \{\pm1\}^n. 
\end{equation}   
\begin{lemma} \label{lem:g} The polynomial in $(\ref{eqn:g})$ is 
expressed as  
\begin{equation} \label{eqn:g-M}
g_m^{(k)}(x) = \sum_{A} \dfrac{(A \mathbf{1})!}{A!} x^{\vec{e} A}, 
\end{equation}
where the sum is taken over all $(k+1) \times n$ quadrilateral matrices 
$A$ of nonnegative integers whose entries sum up to $m$ and 
moreover whose column-sums are all even.     
\end{lemma}
{\it Proof}. Substituting formula (\ref{eqn:h-M}) into definition 
(\ref{eqn:g}) yields 
\begin{equation} \label{eqn:g-M2}
g_m^{(k)}(x) = \dfrac{1}{2^n} \sum_{\varepsilon \in \{\pm 1\}^n} 
\sum_{A} \dfrac{(A \mathbf{1})!}{A!} \varepsilon^{\vec{e} A} x^{\vec{e} A} 
= \dfrac{1}{2^n} \sum_{A} \dfrac{(A \mathbf{1})!}{A!} 
\left( \sum_{\varepsilon \in \{\pm 1\}^n} \varepsilon^{\vec{e} A} \right) 
x^{\vec{e} A},     
\end{equation} 
where the matrix $A$ ranges in the same manner as in formula 
(\ref{eqn:h-M}). 
Put $\vec{e} A = (\nu_1, \dots, \nu_n)$, where $\nu_j$ is 
the $j$-th column-sum of $A$. 
Observe that 
\[
\sum_{\varepsilon \in \{\pm1\}^n} \varepsilon^{\vec{e} A} = \sum_{\varepsilon \in \{\pm1\}^n} 
\varepsilon_1^{\nu_1} \cdots \varepsilon_n^{\nu_n} 
= \left\{ \begin{array}{ll}
2^n \quad & (\mbox{$\nu_j$ is even for any $j = 1,\dots,n$}), \\[2mm]
0   \quad & (\mbox{$\nu_j$ is odd for some $j = 1,\dots,n$}). 
\end{array} \right. 
\]  
So the sum $\sum_A$ in (\ref{eqn:g-M2}) can be restricted to 
those $A$'s whose column-sums are all even.  \hfill $\Box$ 
\par\medskip
For any matrix with even column-sums, its entries must sum up to 
an even number, so that formula (\ref{eqn:g-M}) implies that 
$g_m^{(k)}(x)$ vanishes identically for every $m$ odd. 
Thus from now on $m$ is replaced by $2m$ with $m$ being a positive integer. 
This allows us to put $\vec{e} A = 2 \nu(A)$ with 
$\nu(A) = (\nu_1(A), \dots, \nu_n(A)) \in \mathcal{P}_n(m)$, that is,  
$\nu(A)$ is an ordered $n$-partition of $m$.   
The polynomial $\tau_{2m}^{(k)}(x)$ in formula (\ref{eqn:tau}) is the 
$S_n$-symmetrization of $g_{2m}^{(k)}(x)$, that is, 
\begin{equation} \label{eqn:f-g}
\tau_{2m}^{(k)}(x) = \dfrac{1}{n!} \sum_{\sigma \in S_n} 
g_{2m}^{(k)}(x_{\sigma(1)}, \dots, x_{\sigma(n)}). 
\end{equation}
Putting formula (\ref{eqn:g-M}) with $m$ replaced by $2m$ into 
formula (\ref{eqn:f-g}) yields 
\begin{equation} \label{eqn:tau-M}
\tau_{2m}^{(k)}(x) = \dfrac{1}{n!} \sum_{A \in \mathcal{M}_{n,m}^{(k)}} 
\dfrac{(A \mathbf{1})!}{A!} \sum_{\sigma \in S_n} x_{\sigma(1)}^{2\nu_1(M)} \cdots 
x_{\sigma(n)}^{2\nu_n(A)},   
\end{equation}
where $\mathcal{M}_{n,m}^{(k)}$ is the set of all $(k+1) \times n$ upper 
quadrilateral matrices of nonnegative integers whose entries 
sum up to $2m$ and moreover whose column-sums are all even. 
\par
Let $\zeta_j$ denote a primitive $j$-th root of unity. 
Since 
\[
e_{2j}(\zeta_{2m}, \zeta_{2m}^2, \dots, \zeta_{2m}^m, 
\overbrace{0, \dots, 0}^{n-m}) 
= \left\{ 
\begin{array}{ll}
0 \quad & (j = 1, \dots, m-1), \\[2mm]
(-1)^{m-1} \quad  & (j = m),  
\end{array} \right. 
\]
substituting $x = (\zeta_{2m}, \zeta_{2m}^2, \dots, \zeta_{2m}^m, 
0, \dots, 0)$ into equation (\ref{eqn:tau-s}) yields  
\begin{equation} \label{eqn:c}
c_{n,m}^{(k)} = (-1)^{m-1} 
\tau_{2m}^{(k)}(\zeta_{2m}, \zeta_{2m}^2, \dots, \zeta_{2m}^m, 
0, \dots, 0).   
\end{equation}
For each partition $\nu = (\nu_1,\dots,\nu_n) \in \mathcal{P}_n(m)$, we define       
\begin{equation} \label{eqn:u}
u_{n,m}(\nu) = u_m(\nu_1, \dots, \nu_n) := \sum_{\sigma \in S_{n,m}(\nu)} 
\zeta_m^{\sigma(1)\nu_1+\cdots+\sigma(n)\nu_n},   
\end{equation}
where $S_{n,m}(\nu) := \{\sigma \in S_n \,:\, \mbox{for $i = 1,\dots, n$, 
if $\nu_i \ge 1$ then $\sigma(i) \in \{1,\dots,m \}$ } \}$. 
Since 
\[
x_{\sigma(1)}^{2\nu_1} \cdots x_{\sigma(n)}^{2\nu_n} = \left\{
\begin{array}{ll}
\zeta_m^{\sigma(1)\nu_1+\cdots+\sigma(n)\nu_n} \quad & 
(\sigma \in S_{n,m}(\nu)), \\[2mm] 
0 \quad & (\sigma \not\in S_{n,m}(\nu)),
\end{array}\right. 
\]
at $x = (\zeta_{2m}, \zeta_{2m}^2, \dots, \zeta_{2m}^m, 
0, \dots, 0)$, formulas (\ref{eqn:tau-M}) and (\ref{eqn:c}) yield 
\begin{equation} \label{eqn:c2} 
c_{n,m}^{(k)} = \dfrac{(-1)^{m-1}}{n!} 
\sum_{A \in \mathcal{M}_{n,m}^{(k)}} u_{n,m}(\nu(A)) \dfrac{(A \mathbf{1})!}{A!}.    
\end{equation}
\par
Let $\ell(\nu)$ denote the number of positive entries in an ordered 
partition  $\nu = (\nu_1, \dots, \nu_n) \in \mathcal{P}_n(m)$. 
Note that $\ell(\nu) \in \{1,\dots,m\}$, because $\ell(\nu) \le 
\nu_1 + \cdots + \nu_n = m$.    
\begin{lemma} \label{lem:u-nm}
The function $u_{n,m}(\nu)$ is symmetric, that is, invariant under 
any permutation of $\nu_1,\dots,\nu_n$. 
For any $\nu =(\nu_1,\dots,\nu_n) \in \mathcal{P}_n(m)$ such that  
$\nu_{m+1} = \cdots = \nu_n = 0$, we have   
\begin{equation} \label{eqn:u-nm}
u_{n,m}(\nu) = \dfrac{(n-\ell(\nu))!}{(m-\ell(\nu))!} \, 
u_{m,m}(\nu_1, \dots, \nu_m). 
\end{equation}
\end{lemma}
{\it Proof}. 
For an element $\tau \in S_n$ put $\nu^{\tau} := (\nu_{\tau(1)}, \dots, 
\nu_{\tau(n)})$. 
Then it is easy to see that $S_{n,m}(\nu^{\tau}) = S_{n,m}(\nu) 
\cdot \tau$.  
Using this we show that $u_{n,m}(\nu^{\tau}) = u_{n,m}(\nu)$. 
Indeed, 
\[
\begin{array}{rclcl}
u_{n,m}(\nu^{\tau}) &=& \displaystyle \sum_{\sigma \in S_{n,m}(\nu^{\tau})} 
\zeta_m^{\sigma(1) \nu_{\tau(1)} + \cdots + \sigma(n) \nu_{\tau(n)}}  
&=& \displaystyle \sum_{\sigma \in S_{n,m}(\nu^{\tau})} 
\zeta_m^{(\sigma \cdot \tau^{-1})(1) \cdot \nu_1 + \cdots + 
(\sigma \cdot \tau^{-1})(n) \cdot \nu_n} \\[6mm] 
&=& \displaystyle \sum_{\sigma' \in S_{n,m}(\nu^{\tau}) \cdot \tau^{-1}} 
\zeta_m^{\sigma'(1) \nu_1 + \cdots + \sigma'(n) \nu_n}  
&=& \displaystyle \sum_{\sigma' \in S_{n,m}(\nu)} 
\zeta_m^{\sigma'(1) \nu_1 + \cdots + \sigma'(n) \nu_n} = u_{n,m}(\nu),  
\end{array}
\]
as desired. 
This proves that $u_{n,m}(\nu)$ is a symmetric function of 
$\nu = (\nu_1, \dots, \nu_n)$.  
\par 
We proceed to the second assertion. 
Suppose that $\nu$ is of the form 
$\nu = (\nu_1, \dots, \nu_{r}, 0, \dots, 0)$ with $r := \ell(\nu) \le m$. 
Then $S_{n,m}(\nu) = \{ \sigma \in S_n \,:\, 
\sigma(\{1, \dots, r\}) \subset \{1, \dots, m\} \}$. 
We think of $S_m$ as a subgroup of $S_n$ by setting 
$S_m := \{ \sigma \in S_n \,:\, \sigma(i) = i \,\, \mbox{for} \,\, 
i = m+1, \dots, n\}$. 
Define a map 
\begin{equation} \label{eqn:si-tau} 
S_{n,m}(\nu) \to S_m, \quad \sigma \mapsto \tau \qquad 
\mbox{by} \quad  
\tau(i) := \left\{\begin{array}{ll}
\sigma(i) \quad & (i = 1, \dots, r), \\[2mm]
p(i)   \quad & (i = r+1, \dots, m), \\[2mm]
i      \quad & (i = m+1, \dots, n), 
\end{array}\right.
\end{equation}
where $p$ is the unique bijection $p : \{r+1, \dots, m\} 
\to \{1,\dots,m\} \setminus \sigma(\{1, \dots, r\})$ which 
is ``order-equivalent" to the injection 
$\sigma|_{\{r+1, \dots, m\}}$ in the sense that 
$p(i) < p(j)$ if and only if 
$\sigma(i) < \sigma(j)$ for every $i, j \in \{r+1, \dots, m\}$. 
We claim that the map (\ref{eqn:si-tau}) is 
$\frac{(n-r)!}{(m-r)!}$-to-one. 
Indeed, given any element $\tau \in S_m$, the fiber over 
$\tau$ has a one-to-one correspondence with the set of data 
$(S,q)$:  
\begin{itemize}
\item a subset $S$ of cardinality $m-r$ of $T := 
\{1, \dots, n \} \setminus \tau(\{1, \dots, r\})$, 
\item a bijection $q : \{m+1, \dots, n\} \to T \setminus S$.  
\end{itemize}
It is clear from definition (\ref{eqn:si-tau}) that given a 
data $(S,q)$ there exists a unique element 
$\sigma \in S_{n,m}(\nu)$ such that $\sigma(\{r+1, \dots, m\}) = 
S$ and $\sigma|_{\{m+1, \dots, n\}} = q$.  
Since $\# T = n-r$, there are ${n-r \choose m-r}$ choices of 
$S$, for each of which there are $(n-m)!$ choices of $q$. 
Thus the fiber has a total of ${n-r \choose m-r} \cdot 
(n-m)! = \frac{(n-r)!}{(m-r)!}$ elements. 
Since $\zeta_m^{\sigma(1) 
\nu_1 + \cdots + \sigma(n) \nu_n} = \zeta_m^{\tau(1) \nu_1 + 
\cdots + \tau(m) \nu_m}$, we have    
\[
u_{n,m}(\nu) = \dfrac{(n-r)!}{(m-r)!} 
\sum_{\tau \in S_m} \zeta_m^{\tau(1) \nu_1 + \cdots + \tau(m) \nu_m} 
=  \dfrac{(n-r)!}{(m-r)!} \, u_{m,m}(\nu_1, \dots, \nu_m), 
\]
where $S_m = S_{m,m}(\nu_1, \dots, \nu_m)$ is used in 
the second equality.  \hfill $\Box$ \par\medskip
Formula (\ref{eqn:u-nm}) reduces the calculation of 
$u_{n,m}$ to that of $u_{m,m}$, which we now carry out.    
\begin{lemma} \label{lem:u-nn}
For any partition $\nu = (\nu_1,\dots,\nu_m) \in \mathcal{P}_m(m)$, 
\begin{equation} \label{eqn:u-mm}
u_m(\nu) := u_{m,m}(\nu) = 
m \, (-1)^{\ell(\nu)-1} \, (\ell(\nu)-1)! \, (m-\ell(\nu))!.   
\end{equation}
\end{lemma}
{\it Proof}. 
When $n =m$ the function $u_{n,m}(\nu)$ in definition (\ref{eqn:u}) 
becomes simpler because $S_{m,m}(\nu) = S_m$ for every $\nu 
\in \mathcal{P}_m(m)$.  
The proof is by induction on $\ell(\nu)$. 
When $\ell(\nu) = 1$, we may assume that $\nu$ is of the form 
$\nu = (m,0,\dots,0)$ by the symmetry of $u_m(\nu)$. 
Then definition (\ref{eqn:u}) reads  
\[
u_m(\nu) = \sum_{\sigma \in S_m} \zeta_m^{\sigma(1) m} = 
\sum_{\sigma \in S_m} 1 = m!,  
\]
which verifies formula (\ref{eqn:u-mm}) for $\ell(\nu) = 1$. 
Let $1 \le r < m$ and assume that formula (\ref{eqn:u-mm}) is 
true for every partition $\nu \in \mathcal{P}_m(m)$ with $\ell(\nu) = r$. 
Consider the case $\ell(\nu) = r+1$. 
By the symmetry of $u_m(\nu)$ we may assume that $\nu$ is of the 
form $\nu = (\nu_1,\dots,\nu_{r+1},0,\dots,0)$ with  
$\nu_1, \dots, \nu_{r+1} \ge 1$ and $\nu_1+\cdots+\nu_{r+1} = m$. 
Note that $1 \le \nu_{r+1} < m$.   
Formula (\ref{eqn:u}) now reads 
\[
u_m(\nu) = \sum_{\sigma \in S_m} 
\zeta_m^{\sigma(1)\nu_1+\cdots+\sigma(r)\nu_{r+1}} 
= (m-r-1)! \sum_{(p_1,\dots,p_{r+1})} 
\zeta_m^{p_1 \nu_1+\cdots+p_{r+1} \nu_{r+1}},  
\]
where $(p_1,\dots,p_{r+1})$ ranges over all permutations of distinct 
$r+1$ numbers in $\{1,\dots,m\}$. 
Thus, 
\[
\begin{array}{rcl}
u_m(\nu) &=& (m-r-1)! \displaystyle \sum_{(p_1,\dots,p_r)} 
\zeta_m^{p_1 \nu_1+\cdots+p_r \nu_r} 
\displaystyle \sum_{p_{r+1} \in \{1,\dots,m\} \setminus \{p_1,\dots,p_r \}} 
\zeta_m^{p_{r+1} \nu_{r+1}} \\[6mm]
&=& (m-r-1)! \displaystyle \sum_{(p_1,\dots,p_r)} 
\zeta_m^{p_1 \nu_1+\cdots+p_r \nu_r} \left(
\displaystyle \sum_{l = 1}^m \zeta_m^{l \nu_{r+1}} - 
\displaystyle \sum_{j=1}^r \zeta_m^{p_j \nu_{r+1}} \right). 
\end{array}
\]
Since $1 \le \nu_{r+1} < m$, we have 
$\sum_{l=1}^m \zeta_m^{l \nu_{r+1}} = 0$ and hence 
\[
u_m(\nu) = - (m-r-1)! \displaystyle \sum_{j=1}^r \sum_{(p_1,\dots,p_r)} 
\zeta_m^{p_1 \nu_1^{(j)}+\cdots+p_r \nu_r^{(j)}} = 
- \dfrac{1}{m-r} \sum_{j=1}^r v(\nu^{(j)}),   
\]
where $\nu^{(j)} = (\nu_1^{(j)}, \dots, \nu_m^{(j)}) := 
(\nu_1+ \delta_{1j} \nu_{r+1}, \, \dots, \, 
\nu_1+ \delta_{rj} \nu_{r+1}, \, 0, \, \dots, \, 0)$ with 
$\delta_{ij}$ Kronecker's symbol.  
Note that for each $j = 1,\dots, r$, we have  
$\nu^{(j)} \in \mathcal{P}_m(m)$ with $r(\nu^{(j)}) = r$, 
so that the induction hypothesis yields  
$u_m(\nu^{(j)}) = m \, (-1)^{r-1} \, (r-1)! \, (m-r)!$ 
for $j = 1, \dots, r$.  
Therefore, 
\[
u_m(\nu) = - \dfrac{1}{m-r} \, r \,  
m \, (-1)^{r-1} \, (r-1)! \, (m-r)! = 
m \, (-1)^r \, r! \, (m-r-1)!, 
\]
which means that formula (\ref{eqn:u-mm}) 
is true for $\ell(\nu) = r+1$.  
The induction is complete. \hfill $\Box$
\par\medskip
A column of a matrix is said to be {\sl nontrivial} if 
it has at least one nonzero entry.  
\begin{proposition} \label{prop:c3}
Let $\ell(A)$ denote the number of nontrivial columns in $A$. 
Then, 
\begin{equation} \label{eqn:c3}
c_{n,m}^{(k)} = \dfrac{(-1)^{m-1} \, m}{n!} \sum_{A \in \mathcal{M}_{n,m}^{(k)}} 
(-1)^{\ell(A)-1} \,  
(\ell(A)-1)! \, (n - \ell(A))! \, \dfrac{(A \mathbf{1})!}{A!}.  
\end{equation}
\end{proposition} 
{\it Proof}. 
First, Lemmas \ref{lem:u-nm} and \ref{lem:u-nn} are put together to 
yield the formula  
\begin{equation} \label{eqn:u-nm2}
u_{n,m}(\nu) = m \, (-1)^{\ell(\nu)-1} \, (\ell(\nu)-1)! \, (n-\ell(\nu))!,    
\end{equation}
for any partition $\nu = (\nu_1,\dots,\nu_n) \in \mathcal{P}_n(m)$.  
Indeed, by the symmetry of $u_{n,m}(\nu)$ we may assume  
$\nu_{m+1} = \cdots = \nu_n = 0$. 
Thus using formula (\ref{eqn:u-mm}) in formula (\ref{eqn:u-nm}) gives 
formula (\ref{eqn:u-nm2}).  
Next, putting formula (\ref{eqn:u-nm2}) with 
$\nu = \nu(A)$ into (\ref{eqn:c2}) yields formula (\ref{eqn:c3}),  
since $\ell(A) = \ell(\nu(A))$.  \hfill $\Box$
\section{Young Diagram Representation} \label{sec:young}
We rewrite formula (\ref{eqn:c3}) in Proposition \ref{prop:c3} as a sum 
over some Young diagrams. 
After several preliminary discussions, the main result of this section 
is stated in Proposition \ref{prop:c5}.  
For each $\nu = (\nu_1, \dots, \nu_n) \in \mathbb{Z}_{\ge0}^n$, let 
$\mathcal{M}_n^{(k)}(\nu)$ be the set of all $(k+1) \times n$ upper quadrilateral 
matrices the $i$-th column of which sums up to $\nu_i$ for 
$i = 1, \dots, n$. 
Motivated by expression (\ref{eqn:c3}), put
\begin{equation} \label{eqn:v-def}
v_{n}^{(k)}(\nu) := 
\sum_{A \in \mathcal{M}_n^{(k)}(\nu)} \dfrac{(A \mathbf{1})!}{A!}.  
\end{equation}
\begin{lemma} \label{lem:v}
For any $\nu = (\nu_1, \dots, \nu_n) \in \mathbb{Z}_{\ge0}^n$, we have  
\begin{equation} \label{eqn:v}
v_n^{(k)}(\nu) = 
\dfrac{(\nu_1 + \cdots + \nu_n + k)!}{ 
\prod_{j=1}^k(\nu_1+ \cdots + \nu_j+j) \cdot \prod_{j=1}^n \nu_j!}. 
\end{equation}
\end{lemma}
{\it Proof}. The proof is by induction on $k$. 
For $k = 0$ there is nothing to prove. 
Suppose that formula (\ref{eqn:v}) is true for $k-1$.  
We write $\nu = \vec{\nu}$ to emphasize that $\nu$ is a 
row vector. 
Put 
\[
\psi_n^{(k)}(\vec{a}_1, \dots, \vec{a}_{k+1}) := 
\dfrac{(A \mathbf{1})!}{A!} 
\quad \mbox{for} \quad 
A = \begin{pmatrix} 
\vec{a}_1 \\ \vdots \\ \vec{a}_{k+1} 
\end{pmatrix} 
\quad \mbox{with} \quad 
\vec{a}_i = (\overbrace{0, \dots, 0}^{i-1}, a_{ii}, \dots, a_{in}).   
\]
Here we also write $\vec{a}_{k+1} = \vec{a} = (0, \dots, 0, a_{k+1}, 
\dots, a_n)$ to simplify the notation.  
Observe that $\psi_n^{(k)}(\vec{a}_1, \dots, \vec{a}_k, \vec{a}) = 
\psi_n^{(0)}(\vec{a}) \cdot \psi_n^{(k-1)}(\vec{a}_1, \dots, \vec{a}_k)$.  
Using this we have    
\[
\begin{array}{rcl}
v_n^{(k)}(\vec{\nu}) 
&=& \displaystyle \sum_{\vec{a}_1+\cdots+\vec{a}_k+ \vec{a} =\vec{\nu}} 
\psi_n^{(k)}(\vec{a}_1,\dots, \vec{a}_k, \vec{a}) \\[7mm]  
&=& \displaystyle \sum_{\vec{a} \le \vec{\nu}} 
\psi_n^{(0)}(\vec{a}) 
\displaystyle \sum_{\vec{a}_1+\cdots+\vec{a}_k = \vec{\nu}- \vec{a}} 
\psi_n^{(k-1)}(\vec{a}_1,\dots, \vec{a}_k) 
= \displaystyle \sum_{\vec{a} \le \vec{\nu}} 
\psi_n^{(0)}(\vec{a}) \cdot v_n^{(k-1)}(\vec{\nu}- \vec{a}),  
\end{array}
\]
where $\vec{a} \le \vec{\nu}$ means that  
$\vec{\nu} - \vec{a} \in \mathbb{Z}_{\ge0}^n$. 
Put $\mu_j := \nu_1+ \cdots+\nu_j$, $\bar{\mu}_j := \nu_{j+1} + \cdots + \nu_n$ 
and $b := a_{k+1} + \cdots + a_{n}$.  
Since $a_j = 0$ for $j = 1, \dots, k$, the induction hypothesis 
yields     
\[
v_n^{(k-1)}(\vec{\nu}-\vec{a}) = 
\dfrac{1}{\prod_{j=1}^{k-1} (\mu_j+j) \cdot \prod_{j=1}^k \nu_j!}  
\cdot    
\dfrac{(\mu_n+k-1-b)!}{\prod_{j=k+1}^n (\nu_j-a_j)!}.  
\]
Substituting this into the previous formula and after some 
manipulations we have 
\[
\begin{array}{rcl}
v_n^{(k)}(\vec{\nu}) &=&   
\dfrac{1}{\prod_{j=1}^{k-1} (\mu_j+j) \cdot \prod_{j=1}^n \nu_j!} 
\displaystyle \sum_{b=0}^{\bar{\mu}_k} b! \, (\mu_n+k-1-b)! 
\sum_{\scriptsize \begin{array}{c} \vec{a} \le \vec{\nu} \\[-0.5mm] 
a_{k+1}+\cdots+a_n = b \end{array}} 
\prod_{j=k+1}^n {\nu_j \choose a_j} \\[6mm] 
&=& 
\dfrac{1}{\prod_{j=1}^{k-1} (\mu_j+j) \cdot \prod_{j=1}^n \nu_j!} 
\displaystyle \sum_{b=0}^{\bar{\mu}_k} b! \, (\mu_n+k-1-b)! \, {\bar{\mu}_k \choose b} \\[6mm] 
&=& 
\dfrac{\bar{\mu}_k! \, (\mu_k+k-1)!}{\prod_{j=1}^{k-1} (\mu_j+j) \cdot \prod_{j=1}^n \nu_j!} 
\displaystyle \sum_{b=0}^{\bar{\mu}_k} {\mu_n+k-1-b \choose \mu_k+k-1} \\[7mm]  
&=& 
\dfrac{\bar{\mu}_k! \, (\mu_k+k-1)!}{\prod_{j=1}^{k-1} (\mu_j+j) \cdot \prod_{j=1}^n \nu_j!} 
\displaystyle {\mu_n+k \choose \mu_k+k} = 
\dfrac{(\mu_n+k)!}{\prod_{j=1}^k (\mu_j+j) \cdot \prod_{j=1}^n \nu_j!},   
\end{array}
\]       
where the following general formulas are used to obtain 
the second and fourth  equalities. 
\[
\sum_{\scriptsize \begin{array}{c} 0 \le b_i \le a_i \\ b_1+\cdots+b_k = b 
\end{array}} \prod_{i=1}^k {a_i \choose b_i} = {a_1 + \cdots + a_k \choose b}, 
\qquad 
\sum_{i=b}^a {i \choose b} = {a+1 \choose b+1}. 
\]
Therefore formula (\ref{eqn:v}) remains true for $k$ and  
the induction is complete. \hfill $\Box$
\begin{lemma} \label{lem:c4} 
Formula $(\ref{eqn:c3})$ in Proposition $\ref{prop:c3}$ is 
rewritten as    
\begin{equation} \label{eqn:c4}
c_{n,m}^{(k)} = \dfrac{(-1)^{m-1} m \cdot (2m+k)!}{n!} 
\sum_{\nu \in \mathcal{P}_n(m)} (-1)^{\ell(\nu)-1} \,  
(\ell(\nu)-1)! \, (n - \ell(\nu))! \, \bar{v}_n^{(k)}(\nu),  
\end{equation}
where
\begin{equation} \label{eqn:v-bar}
\bar{v}_n^{(k)}(\nu) := 
\dfrac{1}{\prod_{j=1}^k(2\nu_1+ \cdots + 2\nu_j+j) \cdot 
\prod_{j=1}^n (2\nu_j)!}.   
\end{equation}
\end{lemma}
{\it Proof}. 
Since there exists a direct sum decomposition 
\[
\mathcal{M}_{n,m}^{(k)} = \coprod_{\nu \in \mathcal{P}_n(m)} \mathcal{M}_n^{(k)}(2 \nu),  
\]
and $\ell(A) = \ell(\nu)$ for $A \in \mathcal{M}_n^{(k)}(2\nu)$ 
with $\nu \in \mathcal{P}_n(m)$, formulas (\ref{eqn:c3}) and definition 
(\ref{eqn:v-def}) lead to   
\[
\begin{array}{rcl}
c_{n,m}^{(k)} &=& \dfrac{(-1)^{m-1} \, m}{n!} \displaystyle 
\sum_{\nu \in \mathcal{P}_n(m)} (-1)^{\ell(\nu)-1} \,  
(\ell(\nu)-1)! \, (n - \ell(\nu))! \sum_{A \in \mathcal{M}_n^{(k)}(2\nu)} 
\dfrac{(A \mathbf{1})!}{A!} \\[6mm]
&=& \dfrac{(-1)^{m-1} \, m}{n!} \displaystyle \sum_{\nu \in \mathcal{P}_n(m)} 
(-1)^{\ell(\nu)-1} \,  (\ell(\nu)-1)! \, (n - \ell(\nu))! \, 
v_n^{(k)}(2\nu).  
\end{array} 
\]
Use formula (\ref{eqn:v}) with $\nu$ replaced by $2\nu$ and 
factor the term $(2\nu_1+\cdots+2\nu_n+k)! = (2m+k)!$,  
which is constant for $\nu \in \mathcal{P}_n(m)$, out of the summation.  
Then we obtain formula (\ref{eqn:c4}). \hfill $\Box$ \par\medskip 
Formula (\ref{eqn:c4}) comes up as a sum over ordered partitions. 
The next task is to recast it to a sum over unordered partitions, 
that is, over Young diagrams. 
Let $\mathcal{Y}_j$ be the set of all weakly decreasing sequences 
$\lambda = (\lambda_1 \ge \dots \ge \lambda_j)$ of nonnegative integers.   
The sum $|\lambda| := \lambda_1+\cdots+\lambda_j$ is called the weight of $\lambda$.  
Note that $\lambda$ represents an unordered partition of $|\lambda|$ 
by $j$ nonnegative integers. 
The number of positive entries in $\lambda$, denoted $\ell(\lambda)$, is 
called the length of $\lambda$.   
An element $\lambda \in \mathcal{Y}_j$ defines a Young diagram of weight $|\lambda|$ 
and of length $\ell(\lambda) \le j$ in the usual manner (see e.g. 
Macdonald \cite{Macdonald}).  
An element $\lambda \in \mathcal{Y}_j$ is also written 
$\lambda = \langle 0^{r_0} \, 1^{r_1} \, 2^{r_2} \, \cdots \rangle$ when the 
number $i$ occurs exactly $r_i$ times in $\lambda$  for 
each $i = 0, 1,2, \dots$, where the term $i^{r_i}$ may be omitted 
if $r_i = 0$. 
Note that $|\lambda| = r_1 + 2 r_2 + 3 r_3 + \cdots$,  
$\ell(\lambda) = r_1 + r_2 + r_3 + \cdots$, and $j = r_0 + \ell(\lambda)$.  
Given an element $\lambda \in \mathcal{Y}_j$ let $\mathcal{P}_j(\lambda)$ denote the fiber 
over $\lambda$ of the order-forgetful mapping   
\[
\mathbb{Z}_{\ge0}^j \to \mathcal{Y}_j, \quad \nu = (\nu_1, \dots, \nu_j) \mapsto 
\lambda = \langle 0^{s_0} \, 1^{s_1} \, 2^{s_2} \, \cdots \rangle, 
\]
where $i$ occurs $s_i$ times in $\nu$ for each $i = 0, 1,2, \dots$. 
Motivated by expression (\ref{eqn:v-bar}), consider   
\begin{equation} \label{eqn:w-def}
w_k(\mu) := \sum_{\nu \in \mathcal{P}_k(\mu)} 
\dfrac{1}{\prod_{j=1}^k (2\nu_1 + \cdots + 2\nu_j + j) \cdot 
\prod_{j=1}^k (2\nu_j)!}     
\end{equation}
for $\mu \in \mathcal{Y}_k$, where the denominator of the summand differs 
from that of formula (\ref{eqn:v-bar}) by the factor 
$\prod_{j=1}^k (2\nu_j)!$ in place of $\prod_{j=1}^n (2\nu_j)!$. 
This function is evaluated in the following manner. 
\begin{lemma} \label{lem:w}
For each $\mu = \langle 0^{s_0} \, 1^{s_1} \, 2^{s_2} \cdots \rangle \in 
\mathcal{Y}_k$, 
we have   
\begin{equation} \label{eqn:w}
w_k(\mu) = \dfrac{1}{\prod_{j \ge 0} s_j!\prod_{j \ge0}((2j+1)!)^{s_j}}. 
\end{equation}  
\end{lemma}
{\it Proof}. The proof is by induction on $k$. 
For $k = 1$ formula (\ref{eqn:w}) holds trivially. 
Suppose that $k \ge 2$ and formula (\ref{eqn:w}) is true for $k-1$. 
Let $\mu = \langle j^{s_j} \,|\, j \in J \rangle \in \mathcal{Y}_k$ with 
$J := \{ j \in \mathbb{Z}_{\ge0}\,:\, s_j \neq 0\}$, where $\sum_{j \in J} j s_j 
= |\mu|$ and $\sum_{j \in J} s_j = k$.  
For each $i \in J$, put $\mu^{(i)} := \langle j^{s_j - \delta_{ij}} 
\,|\, j \in J \rangle$, where $\delta_{ij}$ is Kronecker's symbol. 
Since $\mu^{(i)} \in \mathcal{Y}_{k-1}$, the induction hypothesis implies 
that for each $i \in J$, 
\begin{equation} \label{eqn:w-ind}
w_{k-1}(\mu^{(i)}) = \dfrac{1}{\prod_{j \in J} 
(s_j-\delta_{ij})! \prod_{j \in J} ((2j+1)!)^{s_j-\delta_{ij}} } 
= \dfrac{s_i (2i+1)!}{\prod_{j \in J} 
s_j! \prod_{j \in J} ((2j+1)!)^{s_j}}.    
\end{equation}
Observing that there exists a direct sum decomposition 
\[
\mathcal{P}_k(\mu) = \coprod_{i \in J} \{\nu = (\nu_1, \dots, \nu_k) 
\in \mathcal{P}_k(\mu) \,:\, \nu_k = i\} 
= \coprod_{i \in J} \{\nu = (\nu^{(i)}, i) \,:\, \nu^{(i)} 
\in \mathcal{P}_{k-1}(\mu^{(i)}) \}, 
\]
and noticing that $2 \nu_1+\cdots+2\nu_k+k = 2|\mu|+k$, we have 
\begin{eqnarray*}
w_k(\mu) &=& \displaystyle \sum_{i \in J} 
\dfrac{1}{(2|\mu|+k) \cdot (2i)!} \sum_{\nu^{(i)} \in \mathcal{P}_{k-1}(\mu^{(i)})} 
\dfrac{1}{\prod_{j=1}^{k-1}(2\nu_1^{(i)}+\cdots+2\nu_j^{(i)}+j) 
\prod_{j=1}^{k-1} (2\nu_j^{(i)})!} 
\\
&=& \displaystyle \sum_{i \in J} 
\dfrac{1}{(2|\mu|+k) \cdot (2i)!} \, w_{k-1}(\mu^{(i)}) 
\\
&=& \displaystyle \sum_{i \in J} 
\dfrac{1}{(2|\mu|+k) \cdot (2i)!} \cdot 
\dfrac{s_i (2i+1)!}{\prod_{j \in J} 
s_j! \prod_{j \in J} ((2j+1)!)^{s_j}} \hspace{10mm} 
(\mbox{by formula (\ref{eqn:w-ind})}) 
\\
&=& \dfrac{1}{(2|\mu|+k)} \cdot  
\dfrac{\sum_{i \in J} s_i (2i+1)}{\prod_{j \in J} 
s_j! \prod_{j \in J} ((2j+1)!)^{s_j}} 
\\
&=& \dfrac{1}{\prod_{j \in J} 
s_j! \prod_{j \in J} ((2j+1)!)^{s_j}} 
\hspace{10mm}  
(\mbox{by $\sum_{i \in J} s_i (2i+1) = 2 |\mu|+k$}).  
\end{eqnarray*}
This shows that formula (\ref{eqn:w}) remains true for $k$ and the 
induction is complete. \hfill $\Box$ \par\medskip
Let $\mathcal{Y}_n(m)$ be the set of all unordered $n$-partitions of $m$ and  
put $\mathcal{Y}(m) := \mathcal{Y}_m(m)$. 
\begin{proposition} \label{prop:c5}
The generating polynomial $G_{n,m}(t)$ in definition 
$(\ref{eqn:Gnm})$ is expressed as 
\begin{equation} \label{eqn:c5}
(-1)^{m-1} \dfrac{G_{n,m}(t)}{(t+1)^n} = m \hspace{-2mm} 
\sum_{\lambda \in \mathcal{Y}(m)} (-1)^{r_1+\cdots+r_m-1} \, 
\dfrac{(r_1+\cdots+r_m-1)!}{r_1!\cdots r_m!} 
\, T_1^{r_1} \cdots T_m^{r_m}   
\end{equation}
for any $n \ge m \ge 1$, where $\lambda = \langle 0^{r_0} \, 1^{r_1} \, \cdots \, 
m^{r_m}\rangle$ and $T_j$ is defined by 
\begin{equation} \label{eqn:T}
T_j := \dfrac{1}{(2j+1)!} \cdot \dfrac{(2j+1)t+1}{t+1} 
\qquad (j = 1,\dots,m). 
\end{equation} 
In particular the rational function $(t+1)^{-n} G_{n,m}(t)$ is 
independent of $n$.   
\end{proposition}
{\it Proof}. 
For each $\lambda \in \mathcal{Y}_n(m)$ let $\mathcal{Y}_k(\lambda)$ be the set of all Young 
subdiagrams $\mu \in \mathcal{Y}_k$ of $\lambda$.  
For each $\mu \in \mathcal{Y}_k(\lambda)$ let $\mathcal{P}_n^{(k)}(\lambda,\mu)$ be the set 
of all $\nu = (\nu_1,\dots,\nu_n) \in \mathcal{P}_n(\lambda)$ such that the cut-off 
to the first $k$ components $\nu' := (\nu_1,\dots,\nu_k)$ belongs to 
$\mathcal{P}_k(\mu)$.  
Let $\lambda = \langle 0^{r_0} \, 1^{r_1} \, 2^{r_2} \, \cdots \rangle \in \mathcal{Y}_n(m)$ 
and $\mu = \langle 0^{s_0} \, 1^{s_1} \, 2^{s_2} \, \cdots 
\rangle \in \mathcal{Y}_k(\lambda)$. 
Note that $r_j = s_j = 0$ for $j > m$. 
Taking $\mu$ away from $\lambda$ induces a skew-diagram 
$\lambda/\mu = \langle 0^{r_0-s_0} \, 1^{r_1-s_1} \, 2^{r_2-s_2} \, \cdots \rangle$ 
with $r_j-s_j \ge 0$ for $j = 0,1,2,\dots$.  
Since $\# \mathcal{P}_{n-k}(\lambda/\mu) = (n-k)!/\prod_{j \ge 0}(r_j-s_j)!$ and 
$\prod_{i=k+1}^n (2\nu_i)! = \prod_{j \ge 0} ((2j)!)^{r_j-s_j}$, we have    
\begin{eqnarray*}
\displaystyle\sum_{\nu \in \mathcal{P}_n^{(k)}(\lambda,\mu)} \hspace{-2mm} \bar{v}_n^{(k)}(\nu) 
&=&  
\displaystyle \sum_{\nu \in \mathcal{P}_k(\mu)} \dfrac{1}{\displaystyle \prod_{1 \le j\le k} \!
(2\nu_1+\cdots+2\nu_j+j)\hspace{-2mm} 
\prod_{1\le j\le k} \!(2\nu_j)!} \cdot  
\dfrac{(n-k)!}{\displaystyle\prod_{j\ge0}(r_j-s_j)! 
\prod_{j\ge 0}((2j)!)^{r_j-s_j}} 
\\
&=& w_k(\mu) \, \dfrac{(n-k)!}{\prod_{j\ge0}(r_j-s_j)! 
\prod_{j\ge 0}((2j)!)^{r_j-s_j}} \qquad (\mbox{by (\ref{eqn:w-def})}) 
\\
&=& \dfrac{1}{\prod_{j\ge 0} s_j! \prod_{j\ge 0}((2j+1)!)^{s_j}} 
\cdot \dfrac{(n-k)!}{\prod_{j\ge0}(r_j-s_j)! 
\prod_{j\ge 0}((2j)!)^{r_j-s_j}} \quad (\mbox{by (\ref{eqn:w})}) 
\\
&=& \dfrac{(n-k)!}{\prod_{j \ge 0} r_j!} \displaystyle \prod_{j\ge 0} 
{r_j \choose s_j} \left(\dfrac{1}{(2j+1)!} \right)^{s_j} 
\left(\dfrac{1}{(2j)!} \right)^{r_j-s_j}
\end{eqnarray*} 
Since there is a direct sum decomposition $\mathcal{P}_n(\lambda) = 
\coprod_{\mu \in \mathcal{Y}_k(\lambda)} \mathcal{P}_n^{(k)}(\lambda,\mu)$, 
we have  
\begin{eqnarray}
\displaystyle \sum_{k=0}^n \dfrac{t^{n-k}}{(n-k)!} \hspace{-1mm} 
\sum_{\nu \in \mathcal{P}_n(\lambda)} \hspace{-2mm} \bar{v}_n^{(k)}(\nu) 
&=&  
\displaystyle \sum_{k=0}^n \dfrac{t^{n-k}}{(n-k)!} \hspace{-1mm}
\displaystyle \sum_{\mu \in \mathcal{Y}_k(\lambda)} \, \sum_{\nu \in \mathcal{P}_n^{(k)}(\lambda,\mu)} 
\hspace{-2mm} \bar{v}_n^{(k)}(\nu) \nonumber
\\ 
&=& \displaystyle \sum_{k=0}^n \dfrac{t^{n-k}}{(n-k)!} 
\hspace{-5mm} \displaystyle \sum_{\tiny \begin{array}{c} 0 \le s_j \le r_j \\ 
s_0+s_1+\cdots=k \end{array}} \hspace{-5mm}
\dfrac{(n-k)!}{\prod_{j \ge 0} r_j!} \displaystyle \prod_{j\ge 0} 
{r_j \choose s_j} \!\! \left(\dfrac{1}{(2j+1)!} \right)^{s_j} \!\! 
\left(\dfrac{1}{(2j)!} \right)^{r_j-s_j} \nonumber
\\ 
&=& \dfrac{1}{\prod_{j=0}^m r_j!} 
\displaystyle \sum_{s_0=0}^{r_0} \cdots \sum_{s_m=0}^{r_m} \displaystyle 
\prod_{j=0}^m \displaystyle {r_j \choose s_j} 
\!\! \left(\dfrac{1}{(2j+1)!} \right)^{s_j} \!\! 
\left(\dfrac{t}{(2j)!} \right)^{r_j-s_j} \nonumber
\\ 
&=& \dfrac{1}{\prod_{j=0}^m r_j!} 
\displaystyle \prod_{j=0}^m 
\left(\dfrac{1}{(2j+1)!} + \dfrac{t}{(2j)!} \right)^{r_j} \nonumber
\\
&=& \dfrac{(1+t)^{r_0}}{r_0! \prod_{j=1}^m r_j!} 
\displaystyle \prod_{j=1}^m (t+1)^{r_j} T_j^{r_j} \nonumber 
\\
&=& \dfrac{(1+t)^{n}}{r_0! \prod_{j=1}^m r_j!} 
\displaystyle \prod_{j=1}^m T_j^{r_j} \label{eqn:T2}
\end{eqnarray}
where $n-k = \sum_{j=0}^m(r_j-s_j)$ and 
$n = \sum_{j=0}^m r_j$ are used in the 
third and final equalities. 
\par
On the other hand, in view of $\mathcal{P}_n(m) = \coprod_{\lambda \in \mathcal{Y}_n(m)} 
\mathcal{P}_n(\lambda), $ formula (\ref{eqn:c4}) yields
\[
(-1)^{m-1} \dfrac{n! \, c_{n,m}^{(k)}}{(2m+k)!} 
= m \sum_{\lambda \in \mathcal{Y}_n(m)} (-1)^{\ell(\lambda)-1} \, 
(\ell(\lambda)-1)! \, 
(n-\ell(\lambda))! \sum_{\nu \in \mathcal{P}_n(\lambda)} \bar{v}_n^{(k)}(\nu).    
\]
Thus, taking $\ell(\lambda) = r_1+\cdots+r_m$ and $n-\ell(\lambda) = r_0$ 
into account, we have 
\begin{eqnarray*} 
(-1)^{m-1} \dfrac{G_{n,m}(t)}{(t+1)^n} 
&=& \dfrac{(-1)^{m-1}}{(t+1)^n} \sum_{k=0}^n 
\dfrac{n! \, c_{n,m}^{(k)}}{(n-k)! \, (2m+k)!} \, t^{n-k} 
\qquad (\mbox{by (\ref{eqn:Gnm})}) \\
&=& \dfrac{m}{(t+1)^n} \hspace{-2mm} \sum_{\lambda \in \mathcal{Y}_n(m)}  
\hspace{-1mm} (-1)^{\ell(\lambda)-1} \, (\ell(\lambda)-1)! \, r_0! \, 
\sum_{k=0}^n \dfrac{t^{n-k}}{(n-k)!} \hspace{-1mm} 
\sum_{\nu \in \mathcal{P}_n(\lambda)} \bar{v}_n^{(k)}(\lambda) \\
&=& m \hspace{-2mm} \sum_{\lambda \in \mathcal{Y}_n(m)} 
\hspace{-1mm} (-1)^{r_1+\cdots+r_m-1} \, 
\dfrac{(r_1+\cdots+r_m-1)!}{r_1! \cdots r_m!} \, 
T_1^{r_1} \cdots T_m^{r_m} \quad (\mbox{by (\ref{eqn:T2})}),     
\end{eqnarray*}
where the sum may be taken over $\mathcal{Y}(m)$, because when $n \ge m$ 
any $\lambda = (\lambda_1 \ge \cdots \ge \lambda_n) \in \mathcal{Y}_n(m)$ is of length 
at most $m$, that is, $\lambda_j = 0$ for 
any $j > m$, and hence $\lambda$ can be identified with 
$\lambda' := (\lambda_1\ge \cdots \ge \lambda_m) \in \mathcal{Y}(m)$.  
This proves formula (\ref{eqn:c5}). 
As the right-hand side of formula (\ref{eqn:c5}) depends only 
on $m$, the rational function $(t+1)^{-n} G_{n,m}(t)$ is 
independent of $n$. \hfill $\Box$
\section{Generating Functions and Bernoulli Numbers} 
\label{sec:generating}
We are now in a position to establish Theorem \ref{thm:generating} and 
Corollary \ref{cor:recursion}. 
\par\medskip\noindent
{\it Proofs of Theorem $\ref{thm:generating}$ and 
Corollary $\ref{cor:recursion}$}.  
Formula (\ref{eqn:Gnm-Gmm}) is an immediate consequence of the 
last assertion in Proposition \ref{prop:c5} that 
$(t+1)^{-n}G_{n,m}(t)$ is independent of $n$. 
The proof of formula (\ref{eqn:Gmm}) is based on the 
following general fact on generating series: if we put 
\[
\beta_m := m \sum_{\lambda \in \mathcal{Y}(m)} (-1)^{r_1+\cdots+r_m-1} 
\dfrac{(r_1+\cdots+r_m-1)!}{r_1!\cdots r_m!} \, 
\alpha_1^{r_1} \cdots \alpha_m^{r_m} \qquad (m = 1,2,3,\dots),    
\]
where $\lambda = \langle 0^{r_0} \, 1^{r_1} \, \cdots m^{r_m} \rangle$, then 
there exists a formal power series expansion 
\begin{equation} \label{eqn:log}
\log\left(1 + \sum_{m=1}^{\infty} \alpha_m z^{2m}\right) = 
\sum_{m=1}^{\infty} \dfrac{\beta_m}{m} \, z^{2m}.     
\end{equation}
We apply this formula to the situation of 
Proposition \ref{prop:c5}, where $\alpha_j = T_j$ in formula (\ref{eqn:T}) 
and 
\begin{equation} \label{eqn:beta}
\beta_m = (-1)^{m-1} (t+1)^{-m} G_m(t) 
\end{equation}  
in formula (\ref{eqn:c5}) with $n = m$. 
We now find  
\begin{eqnarray}
1+\sum_{m=1}^{\infty} T_m z^{2m} 
&=& 1 + \dfrac{t}{t+1} \sum_{m=1}^{\infty}\dfrac{z^{2m}}{(2m)!} 
+ \dfrac{1}{t+1} \sum_{m=1}^{\infty} \dfrac{z^{2m}}{(2m+1)!}  
\nonumber \\
&=& \dfrac{t}{t+1} \sum_{m=0}^{\infty}\dfrac{z^{2m}}{(2m)!} 
+ \dfrac{1}{t+1} \sum_{m=0}^{\infty} \dfrac{z^{2m}}{(2m+1)!} 
\nonumber \\
&=& \dfrac{1}{t+1} \left(t \cosh z + \frac{\sinh z}{z} \right). 
\label{eqn:log2}
\end{eqnarray}
Substitute this into formula (\ref{eqn:log}) and apply the 
differential operator $\frac{z}{2}\frac{\partial}{\partial z}$ 
to the resulting equation.      
Then after some calculations we get formula (\ref{eqn:Gmm}) and 
thus establish Theorem \ref{thm:generating}.  
Corollary \ref{cor:recursion} then follows easily from 
this theorem in the manner mentioned in the Introduction.  
\hfill $\Box$ \par\medskip
We present some amplifications of Theorem \ref{thm:generating} and 
Corollary \ref{cor:recursion}. 
For the extremal cases of $k = 0$, $n-1$, $n$, the coefficients 
$c_{n,m}^{(k)}$ can be written explicitly in terms of 
Bernoulli numbers. 
\begin{lemma} \label{lem:extrema}
For $k = 0$, $n-1$, $n$, the coefficients $c_{n,m}^{(k)}$ are directly 
tied to $b_m$ by  
\begin{equation} \label{eqn:extrema}
c_{n,m}^{(0)} = (2m)! \, (2^{2m}-1) \, b_m 
\qquad 
c_{n,m}^{(n-1)} = c_{n,m}^{(n)} = \dfrac{(n+2m)!}{n!} \, b_m, 
\qquad (n \ge m \ge 1). 
\end{equation} 
\end{lemma}
{\it Proof}. 
Substitute $t = 0$ into definition (\ref{eqn:Gnm}) to get  
$G_{n,m}(0) = \frac{n!}{(n+2m)!} \, c_{n,m}^{(n)}$. 
Similarly put $t = 0$ in formulas (\ref{eqn:Gnm-Gmm}) and 
(\ref{eqn:recursion}) to have $G_{n,m}(0) = G_m(0) = b_m$. 
These together lead to the assertion for $c_{n,m}^{(n)}$ 
in formula (\ref{eqn:extrema}). 
The assertion for $c_{n,m}^{(n-1)}$ then follows from the 
identity $c_{n,m}^{(n-1)} = c_{n,m}^{(n)}$ mentioned in the 
Introduction.     
To prove the assertion for $c_{n,m}^{(0)}$ in formula 
(\ref{eqn:extrema}) we consider the generating polynomial 
$\hat{G}_{n,m}(t) := t^n G_{n,m}(1/t)$ instead of $G_{n,m}(t)$. 
After the change $t \mapsto 1/t$ and multiplication by 
$t^n$, formula (\ref{eqn:Gnm}) gives $\hat{G}_{n,m}(0) = 
c_{n,m}^{(0)}/(2m)!$.  
On the other hand, formula (\ref{eqn:Gnm-Gmm}) yields  
$\hat{G}_{n,m}(t) = (t+1)^{n-m} \hat{G}_m(t)$, where 
$\hat{G}_m(t) := \hat{G}_{m,m}(t)$, while formula 
(\ref{eqn:Gmm}) gives  
\[
\sum_{m=1}^{\infty} (-1)^{m-1} \hat{G}_m(t) 
\left(\dfrac{z^2}{t+1} \right)^m = 
\dfrac{t(z \coth z-1)+z^2}{2(z \coth z + t)},  
\]
which upon putting $t = 0$ reduces to the equality  
\[
\sum_{m=1}^{\infty} (-1)^{m-1} \hat{G}_m(0) z^{2m} 
= \dfrac{z}{2 \coth z} = \dfrac{z}{2} \tanh z. 
\]
Comparing it with the Maclaurin expansion 
$\tanh z = 2 \sum_{m=1}^{\infty} (-1)^{m-1}(2^{2m}-1) 
\, b_m \, z^{2m-1}$, we find 
$\hat{G}_m(0) = (2^{2m}-1) \, b_m$.  
Thus $c_{n,m}^{(0)} = (2m)! \, \hat{G}_{n,m}(0) = 
(2m)! \, \hat{G}_m(0) = (2m)! \, (2^{2m}-1) \, b_m$. 
\hfill $\Box$ \par\medskip 
The first formula in (\ref{eqn:extrema}) is already 
found in \cite{Haeuslein}. 
To deal with the intermediate coefficients $c_{n,m}^{(k)}$ for 
$k = 1,\dots,n-2$, another modification of the generating 
polynomials $G_{n,m}(t)$ is helpful.   
\begin{equation} \label{eqn:Fnm}
F_{n,m}(t) := t^n G_{n,m}\left(\dfrac{1-t}{t}\right) = 
\sum_{k=0}^n \dfrac{n! \, c_{n,m}^{(k)}}{(n-k)! \, (2m+k)!} \,  
t^k (1-t)^{n-k}.    
\end{equation}
\begin{lemma} \label{lem:Fm} 
For $n \ge m \ge 1$, the polynomials $F_{n,m}(t)$ depend only on 
$m$, being independent of $n$.  
They satisfy the differential-difference equation  
\begin{equation} \label{eqn:diff-diff}
2 F_{n,m}(t) + \dfrac{t}{m} F_{n,m}'(t) + 
\dfrac{(1-t)^2}{m-1} F_{n-1,m-1}'(t) = 0 \qquad (n \ge m \ge 2).   
\end{equation}
All the $F_{n,m}(t)$ can be determined inductively by solving 
equation $(\ref{eqn:diff-diff})$ with initial conditions 
\begin{equation} \label{eqn:initial}
F_{n,m}(0) = (2^{2m}-1) b_m, \qquad 
F_{n,1}(t) = \dfrac{1}{2}-\dfrac{t}{3} \qquad (n \ge m \ge 1).  
\end{equation}
\end{lemma} 
{\it Proof}. 
Put $F_m(t) := F_{m,m}(t)$. 
It readily follows from relation (\ref{eqn:Gnm-Gmm}) and definition 
(\ref{eqn:Fnm}) that $F_{n,m}(t) = F_m(t)$ for every $n \ge m$.   
The substitution $t \mapsto \frac{1-t}{t}$ induces the changes 
\[
\beta_m \mapsto (-1)^{m-1} F_m(t), \qquad 
1+\sum_{m=1}^{\infty} T_m z^{2m} 
\mapsto (1-t) \cosh z + t \, \frac{\sinh z}{z} 
\]
in formulas (\ref{eqn:beta}) and (\ref{eqn:log2}) respectively. 
With these changes formula (\ref{eqn:log}) reads
\begin{equation} \label{eqn:log3}
\log\left\{ (1-t) \cosh z + t \, \frac{\sinh z}{z} \right\} 
= \sum_{m=1}^{\infty} \dfrac{(-1)^{m-1}}{m} F_m(t) z^{2m}. 
\end{equation}
Denote the both sides of this equation by $\Phi = \Phi(z,t)$. 
A direct check using the left-hand side of equation (\ref{eqn:log3}) 
tells us that $\Phi$ satisfies the partial differential equation  
\begin{equation} \label{eqn:pde}
z \dfrac{\partial \Phi}{\partial z} + 
\left\{t - (1-t)^2 z^2  \right\} \dfrac{\partial \Phi}{\partial t} 
= (1-t) z^2. 
\end{equation}
Next we look at this equation by means of the right-hand side of 
formula (\ref{eqn:log3}). 
For each $m \ge 2$ the coefficient of $z^{2m}$ in equation 
(\ref{eqn:pde}) being zero gives the differential-difference equation   
\[
2 F_m(t) + \dfrac{t}{m} F_m'(t) + \dfrac{(1-t)^2}{m-1} F_{m-1}'(t) 
= 0 \qquad (m \ge 2),    
\]
which can be expressed as equation (\ref{eqn:diff-diff}), 
because $F_m(t) = F_{n,m}(t)$ and $F_{m-1}(t) = F_{n-1,m-1}(t)$ 
by the first assertion of the lemma. 
The first condition in (\ref{eqn:initial}) is derived  
from formulas (\ref{eqn:Fnm}) and (\ref{eqn:extrema}) as 
$F_{n,m}(0) = c_{n,m}^{(0)}/(2m)! = (2^{2m}-1) b_m$, while  
the second condition follows from $F_{n,1}(t) = F_{1,1}(t)$ 
and the direct calculation of $F_{1,1}(t)$, which is easy.      
\hfill $\Box$ \par\medskip  
Differential-difference equation (\ref{eqn:diff-diff}) can be 
used to derive a recursion formula for $c_{n,m}^{(k)}$ as 
well as to explicitly determine $c_{n,m}^{(k)}$ for $k$ near $0$ or 
$n$ in terms of Bernoulli numbers. 
\begin{proposition} \label{prop:c6} 
For $k = 0$, $1$, the coefficients $c_{n,m}^{(k)}$ are given 
by the first formula in $(\ref{eqn:extrema})$ and 
\begin{equation} \label{eqn:c61}
c_{n,m}^{(1)} = (2m+1)! \left\{(2^{2m}-1) \, b_m - 
\frac{2 m}{n} (2^{2(m+1)}-1) \, b_{m+1} \right\} 
\qquad (n \ge m \ge 1).   
\end{equation}
For $k = n-2$, $n-1$, $n$, the coefficients $c_{n,m}^{(k)}$ 
take a common value which is given by 
\begin{equation} \label{eqn:c62}
c_{n,m}^{(n-2)} = c_{n,m}^{(n-1)} = c_{n,m}^{(n)} = 
\dfrac{(n+2m)!}{n!} \, b_m \qquad (n \ge m \ge 2).  
\end{equation}
Moreover for $1 \le k \le n-2$ and  $2 \le m \le n$ there exists a 
recursion formula 
\begin{equation} \label{eqn:c63}
c_{n,m}^{(k)} - c_{n,m}^{(k-1)} = \dfrac{(n-k)(n-k-1) m}{n(m-1)} 
\left\{(2m+k-1)\, c_{n-1,m-1}^{(k)} - (k+1) \, c_{n-1,m-1}^{(k+1)} 
\right\}. 
\end{equation}
\end{proposition}   
{\it Proof}. Write the left-hand side of equation 
(\ref{eqn:diff-diff}) as $\sum_{k=0}^n \gamma_{n,m}^{(k)} \, 
t^k (1-t)^{n-k}$. 
Since the polynomials $t^k (1-t)^{n-k}$, $k = 0, \dots, n$, form 
a basis of the linear space of polynomials in $t$ of degree at most 
$n$, it follows from equation (\ref{eqn:diff-diff}) that 
$\gamma_{n,m}^{(k)} = 0$ for every $0 \le k \le n$. 
For $k = 0$ we find 
\[
\gamma_{n,m}^{(0)} = 2 \dfrac{c_{n,m}^{(0)}}{(2m)!} + 
\dfrac{n-1}{m-1} 
\left\{\dfrac{c_{n-1,m-1}^{(1)}}{(2m-1)!}-
\dfrac{c_{n-1,m-1}^{(0)}}{(2m-2)!}\right\} = 0 \qquad 
(n \ge m \ge 2),    
\]
where $c_{n,m}^{(0)}$ and $c_{n-1,m-1}^{(0)}$ are already 
known as in the first formula of (\ref{eqn:extrema}). 
Thus $c_{n-1,m-1}^{(1)}$ is also known from this equation. 
Replacing $(n-1,m-1)$ with $(n,m)$ we get formula (\ref{eqn:c61}). 
On the other hand, for $k = n-1$, $n$, some calculations show  
that $\gamma_{n,m}^{(n-1)} = 0$ and $\gamma_{n,m}^{(n)} = 0$ lead to 
$c_{n,m}^{(n-2)} = c_{n,m}^{(n-1)}$ and 
$c_{n,m}^{(n-1)} = c_{n,m}^{(n)}$ respectively, where the 
latter is already pointed out in the Introduction and Lemma 
\ref{lem:extrema}. 
Thus formula (\ref{eqn:c62}) follows from the second formula 
of (\ref{eqn:extrema}). 
Finally some more calculations of $\gamma_{m,n}^{(k)}$ for 
general $k$ imply that equation $\gamma_{n,m}^{(k)} = 0$ with  
$1 \le k \le n-2$ and $2 \le m \le n$ is equivalent to 
recursion formula (\ref{eqn:c63}). \hfill $\Box$ \par\medskip 
Since $c_{n,m}^{(k)}$ is already known for the $k$'s at both ends 
of the interval $0 \le k \le n$ as in formulas (\ref{eqn:extrema}), 
(\ref{eqn:c61}) and (\ref{eqn:c62}), the recursion formula 
(\ref{eqn:c63}) can be used to inductively determine all coefficients 
$c_{n,m}^{(k)}$, where there are three directions in which induction 
works productively. 
\\[2mm]
\qquad (a) \, $c_{n,m}^{(k)} \leftarrow 
c_{n,m}^{(k-1)}, \, c_{n-1,m-1}^{(k)}, \, c_{n-1,m-1}^{(k+1)}$, \quad  
(b) \, $c_{n,m}^{(k-1)} \leftarrow 
c_{n,m}^{(k)}, \, c_{n-1,m-1}^{(k)}, \, c_{n-1,m-1}^{(k+1)}$, \\[2mm] 
\qquad (c) \, $c_{n-1,m-1}^{(k+1)} \leftarrow 
c_{n,m}^{(k)}, \, c_{n,m}^{(k-1)}, \,  c_{n-1,m-1}^{(k)}$ (with 
$(m-1,n-1)$ replaced by $(m,n)$).  
\\[2mm]
For example formula (\ref{eqn:c63}) with $k=n-2$ is used in 
direction (b) to derive         
\[
c_{n,m}^{(n-3)} = \dfrac{1}{n!} 
\left\{(n+2m)! \, b_m - 4m \cdot (n+2m-3)! \, b_{m-1} \right\} 
\qquad (n \ge 3, \, n \ge m \ge 2)   
\]
from formula (\ref{eqn:c62}). 
Similarly formula (\ref{eqn:c63}) with $k=1$ can be applied in 
direction (c) to deduce a closed expression for $c_{m,n}^{(2)}$ 
from formulas (\ref{eqn:extrema}) and (\ref{eqn:c61}), and so on.    
\par
At the end we return to the starting point of 
this article, that is, to polyhedral harmonics. 
With Theorem \ref{thm:harmonic} for the cube case, the determination of 
polyhedral harmonic functions for all skeletons of all regular convex 
polytopes has been completed.  
As a summary we have:       
\begin{theorem} \label{thm:regular} 
Let $P$ be any $n$-dimensional regular convex polytope with center 
at the origin in $\mathbb{R}^n$ and $G$ the symmetry group of $P$.  
Then for any $k = 0, \dots, n$, the linear space $\mathcal{H}_{P(k)}$ is of  
$|G|$-dimensions, where $|G|$ denotes the order of $G$, and as an 
$\mathbb{R}[\partial]$-module $\mathcal{H}_{P(k)}$ is generated by the fundamental 
alternating polynomial $\varDelta_G$ of the reflection group $G$. 
\end{theorem}
\begin{figure}[t] 
\begin{center}
\includegraphics*[width=40mm,height=42mm]{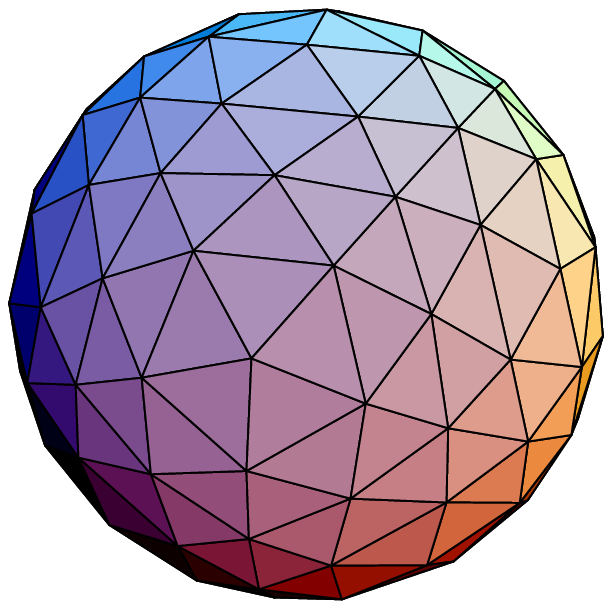}
\includegraphics*[width=40mm,height=42mm]{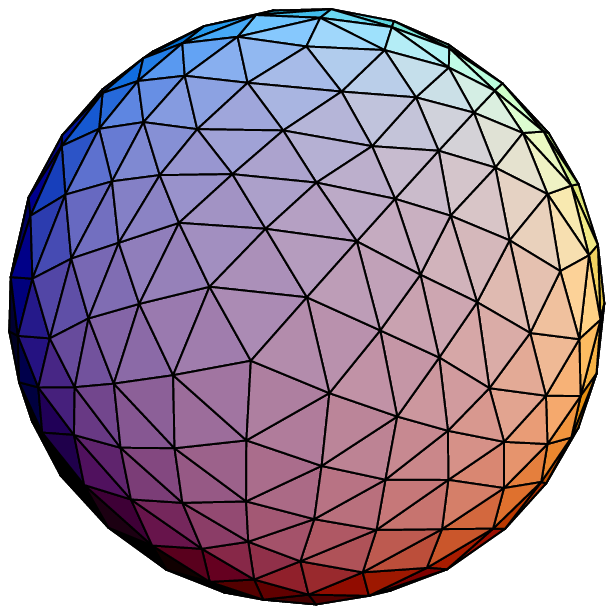}
\includegraphics*[width=40mm,height=42mm]{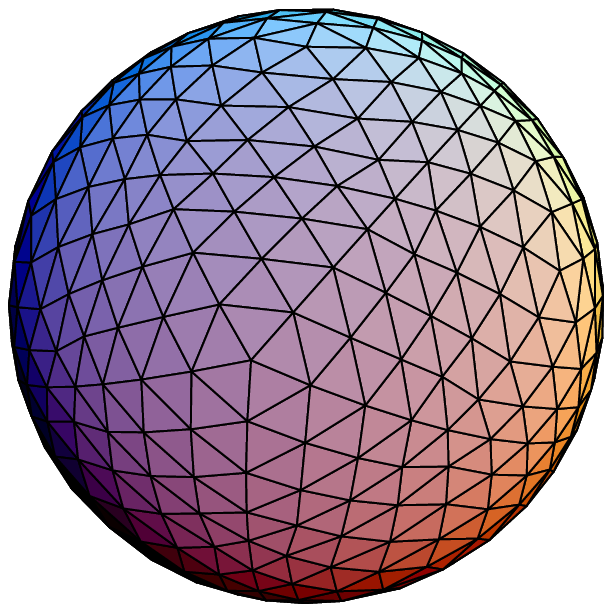} 
\end{center}
\vspace{-7mm} 
\caption{Approximations of the sphere by geodesic domes}
\label{fig:geodesic}
\end{figure}
\par
For the classification of regular convex polytopes we refer to 
Coxeter \cite{Coxeter}. 
Theorem \ref{thm:regular} is proved in article \cite{Iwasaki2} 
for the $n$-dimensional regular simplex and in article \cite{IKM} 
for the exceptional regular polytopes, that is, for the dodecahedron 
and icosahedron in $3$-dimensions and for the $24$-cell, 
$120$-cell and $600$-cell in $4$-dimensions. 
For the $n$-dimensional cross polytope, namely, the analogue 
in $n$-dimensions of the octahedron, there is no detailed written 
proof in the literature, but a proof quite similar to the regular 
$n$-simplex case is feasible.  
This is because each face of an $n$-dimensional cross polytope 
is an $(n-1)$-dimensional regular simplex. 
Finally the $n$-dimensional cube case has been treated in this 
article (Theorem \ref{thm:harmonic}), in which case the proof 
is quite different from those in the other cases. 
Here we should also mention the important studies 
\cite{Flatto1,Flatto2,FW,GR,Haeuslein,Ignatenko} etc. of earlier 
times, which contain partial answers to our questions, referring 
to the survey \cite{Iwasaki3} for a more extensive literature. 
\par
Apart from the regular figures for which symmetry plays a 
dominant role, polyhedral harmonics is largely open, 
for example, for such figures as geodesic domes in 
Figure \ref{fig:geodesic}.  

\end{document}